\documentclass[11pt,hyp,]{amsart}
\usepackage{hyperref}
\hypersetup{nesting=true,debug=true,naturalnames=true}
\usepackage{graphicx,amssymb,upref}

\title{Random walk on a diagonal lattice}
\author{Theo van Uem}  
\address{Amsterdam University of Applied Sciences, Amsterdam, The Netherlands.} 
\email{tjvanuem@gmail.com}  

\subjclass[2020]{60G50,60J05}

\begin{document}
\hbadness=99999


\begin{abstract}
We consider a discrete random walk on a diagonal lattice in two and three dimensions and obtain explicit solutions of
absorption probabilities and probabilities of return in several domains. In three dimensions we consider both the cube and the dodecahedron variant. In two dimensions we obtain explicit formula in case of rotated barriers.
\end{abstract}

\maketitle

\section{Introduction}
Discrete random walks are studied in a number of standard books, see e.g. Spitzer~\cite{SP} and Feller~\cite{FE}.
Polya~\cite{PO} was the first to observe that a walker is certain to return to his starting position in one and two
dimensional symmetric discrete random walks while there exists a positive escape probability in higher
dimensions. McCrea and Whipple~\cite{MC} study simple symmetric random paths in two and three
dimensions, starting in a rectangular lattice on the integers with absorbing barriers on the boundaries.
After taking limits they obtain probabilities of absorption in two and three dimensional lattices.
Bachelor and Henry~\cite{BA1}~\cite{BA2}  use the McCrea-Whipple approach and find the exact solution for random
walks in the triangular lattice with absorbing boundaries and for random walks on finite lattice tubes.
In this paper we study random walks on a diagonal lattice.

\section{Random walk on a diagonal lattice in two dimensions with absorbing boundaries}

\subsection{Rectangular region}
We define an interior $I$ of a rectangular region: $I=\{(p,q)| 1\leq p\leq m, 1\leq q\leq n\}$
The boundary of this region is $B$, which consist of absorbing barriers.
We define $F_{(a,b)}(p,q)$ as the expected number of departures from $(p,q)$ when starting in the interior
source $(a,b)$ on a diagonal lattice. We’ll often use the abbreviation $F(p,q)$.
We study a diagonal lattice, so we have for $I$:
\begin{multline}
\label{eq:one}
F(p,q)=\delta_{a,p} \delta_{b,q} +\\ \frac{1}{4}\{F(p+1,q+1)+F(p+1,q-1)+F(p-1,q+1)+F(p-1,q-1)\}
\end{multline}
and for $B$:

\begin{equation}
\label{eq:two}
F(p,q)=0
\end{equation}
The homogeneous part of the difference equation~\eqref{eq:one} has solutions \(F(p,q)=
Ae^{ip\alpha+q\beta}\), where $\cos\alpha\cosh\beta=1$, so $F(p,q)=C\sin{\alpha p}\sinh{\beta q}.$ \\
We can construct solutions of~\eqref{eq:one} and~\eqref{eq:two}:
\[
F_1(p,q)=\sum_{r=1}^{m} C(r)\sin\frac{pr\pi}{m+1}\sinh q\beta_r\sinh[(n+1-b)\beta_r] \quad (q\leq b)
\]
  \[
F_2(p,q)=\sum_{r=1}^{m} C(r)\sin\frac{pr\pi}{m+1}\sinh b\beta_r\sinh[(n+1-q)\beta_r] \quad (q\geq b)
\]
where \[\cos \frac{r\pi}{m+1}\cosh \beta_r=1\]
We substitute these solutions in~\eqref{eq:one} with $q=b$ and get:
\begin{multline*}
 \sum_{r=1}^{m} C(r)\sin\frac{pr\pi}{m+1}\{\sinh b\beta_r\sinh(n+1-b)\beta_r- \\ \frac{1}{2}\cos\frac{r\pi}{m+1}[\sinh b\beta_r\sinh(n-b)\beta_r+\sinh (b-1)\beta_r\sinh(n+1-b)\beta_r]\}=\delta_{a,p}
 \end{multline*}
 
 Using \(\cos \frac{r\pi}{m+1}\cosh \beta_r=1\) we get after some calculations:

\begin{equation*}
 \sum_{r=1}^{m} C(r)\sin\frac{pr\pi}{m+1}\{ \frac{1}{2}\cos\frac{r\pi}{m+1}\sinh \beta_r\sinh[(n+1)\beta_r]\}=\delta_{a,p}
 \end{equation*}
 
 Using
 \begin{equation*}
 \frac{2}{m+1}\sum_{r=1}^{m} \sin\frac{ar\pi}{m+1}\sin\frac{pr\pi}{m+1}=\delta_{a,p}
 \end{equation*}
 we get:
 \[
F_1(p,q)= \frac{4}{m+1}\sum_{r=1}^{m} \frac{ \sin\frac{ar\pi}{m+1} \sin\frac{pr\pi}{m+1} \sinh q\beta_r\sinh[(n+1-b)\beta_r]}{\tanh \beta_r\sinh[(n+1)\beta_r]} \quad (q\leq b)
\]
  \[
F_2(p,q)= \frac{4}{m+1}\sum_{r=1}^{m} \frac{ \sin\frac{ar\pi}{m+1} \sin\frac{pr\pi}{m+1} \sinh b\beta_r\sinh[(n+1-q)\beta_r]}{\tanh \beta_r\sinh[(n+1)\beta_r]} \quad (q\geq b)
\]
where \[\cos \frac{r\pi}{m+1}\cosh \beta_r=1\]
Remark.
When m is odd we have a problem in $r=\frac{m+1}{2}$. We can change the roles of $p$ and $q$ in our solutions, but when both $m$ and $n$ are odd, our method doesn't work.\\
We obtain absorption probabilities in elements of B by observing (diagonal) neighbors in the interior region. Let $A(p,q)$ be the probability of absorption in $(p,q)\in B$. Then we have for example: $A(m+1,n+1)=\frac{1}{4}{F_2(m,n)},\quad A(m+1,n)=\frac{1}{4}{F_2(m,n-1)},\quad A(m+1,n-1)=\frac{1}{4}[F_2(m,n)+F_2(m,n-2)]$.
\subsection{Semi infinite strip}
By taking $n  \rightarrow \infty$ in the rectangular solution we get:
\[
F_1(p,q)= \frac{4}{m+1}\sum_{r=1}^{m} \frac{ \sin\frac{ar\pi}{m+1} \sin\frac{pr\pi}{m+1} \sinh q\beta_r \exp{(-b\beta_r)}}{\tanh \beta_r} \quad (q\leq b)
\]
 
 \[
F_2(p,q)= \frac{4}{m+1}\sum_{r=1}^{m} \frac{ \sin\frac{ar\pi}{m+1} \sin\frac{pr\pi}{m+1} \sinh b\beta_r \exp{(-q\beta_r)}}{\tanh \beta_r} \quad (q\geq b)
\]
where \[\cos \frac{r\pi}{m+1}\cosh \beta_r=1\]
\subsection{Infinite strip}
By taking $b,q  \rightarrow \infty$ , $q-b=s$ finite in the solution of the semi infinite strip, we get:
\[
F_{(a,0)}(p,s)= \frac{2}{m+1}\sum_{r=1}^{m} \frac{ \sin\frac{ar\pi}{m+1} \sin\frac{pr\pi}{m+1} \exp{(-\lvert s \lvert\beta_r)}}{\tanh \beta_r} 
\]
where \[\cos \frac{r\pi}{m+1}\cosh \beta_r=1\]

\subsection{Infinite Quadrant}
By letting $m,n\rightarrow \infty$ in the block solution, we get the infinite quadrant $p,q>0$.
\[
F_1(p,q)= \frac{8}{\pi}\int_{0}^{\pi} \frac{ \sin{a \lambda} \sin{p \lambda} \sinh{q \mu} \exp{(-b \mu)}}{\tanh \mu} d\lambda \quad (q\leq b)
\]
\[
F_2(p,q)= \frac{8}{\pi}\int_{0}^{\pi} \frac{ \sin{a \lambda} \sin{p \lambda} \sinh{b \mu} \exp{(-b \mu)}}{\tanh \mu} d\lambda \quad (q\geq b)
\]
where 
\begin{equation*}
\cos {\lambda}\cosh {\mu}=1
\end{equation*}
\subsection{Half-plane}

By taking $m  \rightarrow \infty$ ,in the solution of the infinite strip, we get:
\begin{equation}
\label{eq:three}
F_{(a,0)}(p,s)= \frac{2}{\pi}\int_{0}^{\pi} \frac{ \sin{a \lambda} \sin{p \lambda} \exp{(-\lvert s \lvert\mu)}}{\tanh \mu} d\lambda
\end{equation}
where 
\begin{equation}
\label{eq:four}
\cos {\lambda}\cosh {\mu}=1
\end{equation}
We prove this is the unique solution in the half-plane.
First we prove that it is a solution:
If $\lvert s \lvert\geq 1$ then substitute \eqref{eq:three}  in \eqref{eq:one} and use \eqref{eq:four} .
If $s=0$ then we again substitute \eqref{eq:three}  in \eqref{eq:one} and get, using \eqref{eq:four}  with starting point $(a,0)$:
\[
4F(p,0)-F(p+1,1)-F(p+1,-1)-F(p-1,1)-F(p-1,-1)=  
\]
\[
\frac{2}{\pi}\int_{0}^{\pi} \frac{\sin{a\lambda}\left\{ 4\sin{p\lambda}-2 \left[\sin{(p+1)\lambda}+\sin{(p-1)\lambda}  \right] {\rm e}^{-\mu}\right\}}{\tanh \mu}d\lambda=
\]
\[
\frac{2}{\pi}\int_{0}^{\pi} \frac{\sin{a\lambda}\sin{p\lambda}\left[ 4-4\cos\lambda {\rm e}^{-\mu}\right]}{\tanh \mu}d\lambda=\frac{8}{\pi}\int_{0}^{\pi} \sin{a\lambda}\sin{p\lambda}d\lambda=4 \delta_{a,p}
\]
The solution is unique: see Feller [6], (p.362)

\section{Random walk on a diagonal lattice in three dimensions}

\subsection{Block}
The interior is now defined by: $I=\{(p,q,r)| 1\leq p \leq l ,1 \leq q \leq m, 1\leq r\leq n\}$
The boundary of this region is $B$, which consist of absorbing barriers.
We define $F_{(a,b,c)}(p,q,r)$ as the expected number of departures from $(p,q,r)$ when starting in the interior
source $(a,b,c)$ on a diagonal lattice. We’ll often use the abbreviation $F(p,q,r)$.
We study diagonal lattices. In three dimensions this can be realized in two ways: cube and dodecahedron. \\
\textit{We start with the cube model:}

\begin{multline}
\label{eq:five}
F(p,q,r)=\delta_{a,p} \delta_{b,q}  \delta_{c,r} + 
 \frac{1}{8}\{F(p+1,q+1,r+1)+F(p+1,q+1,r-1)+ \\F(p+1,q-1,r+1)+F(p+1,q-1,r-1)+F(p-1,q+1,r+1)+\\F(p-1,q+1,r-1)+F(p-1,q-1,r+1)+F(p-1,q-1,r+1)\}
\end{multline}
and for $B$:
\begin{equation*}
F(p,q,r)=0
\end{equation*}
The homogeneous part of the difference equation~\eqref{eq:five} has solutions
\[(p,q,r)= Ae^{ip\alpha_1+iq\alpha_2+r\beta}\], where $\cos\alpha_1 \cos\alpha_2\cosh\beta=1$, so we have solutions
 \[F(p,q,r)=  C\sin{\alpha_1 p}\sin{\alpha_2 q}\sinh{\beta r}. \]
Analogue to the 2-dimensional case we find  
\begin{multline*}
F_1(p,q,r)= \frac{8}{(l+1)(m+1)}  \\ \sum_{s=1}^{l}\sum_{t=1}^{m} \frac{ \sin\frac{as\pi}{l+1} \sin\frac{ps\pi}{l+1}\sin\frac{bt\pi}{m+1} \sin\frac{qt\pi}{m+1} \sinh r\beta_{st}\sinh[(n+1-c)\beta_{st}]}{\tanh \beta_{st}\sinh[(n+1)\beta_{st}]} \ (r\leq c)
\end{multline*}
\begin{multline*}
F_2(p,q,r)=  \frac{8}{(l+1)(m+1)}  \\ \sum_{s=1}^{l}\sum_{t=1}^{m} \frac{ \sin\frac{as\pi}{l+1} \sin\frac{ps\pi}{l+1}\sin\frac{bt\pi}{m+1} \sin\frac{qt\pi}{m+1} \sinh c\beta_{st}\sinh[(n+1-r)\beta_{st}]}{\tanh \beta_{st}\sinh[(n+1)\beta_{st}]} \ (r\geq c)
\end{multline*}
where \[\cos \frac{s\pi}{l+1}\cos \frac{t\pi}{m+1}\cosh \beta_{st}=1\]
\\
The next model is the \textit{dodecahedron} case:
\begin{multline}
\label{eq:six}
F(p,q,r)=\delta_{a,p} \delta_{b,q}  \delta_{c,r} + \\
 \frac{1}{12}\{F(p+1,q+1,r)+F(p+1,q-1,r)+ F(p-1,q+1,r)+F(p-1,q-1,r)+\\F(p+1,q,r+1)+F(p+1,q,r-1)+F(p-1,q+1,r)+F(p-1,q-1,r)+\\F(p,q+1,r+1)+F(p,q+1,r-1)+F(p,q-1,r+1)+F(p,q-1,r-1)
 \}
 \end{multline}
 
The homogeneous part of the difference equation~\eqref{eq:six} has solutions:\\ \(F(p,q,r)=
Ae^{ip\alpha_1+iq\alpha_2+r\beta}\), where $\cos\alpha_1 \cos\alpha_2+(\cos\alpha_1 +\cos\alpha_2)\cosh\beta=3$.\\ 
The solutions for dodecahedron case are the same as for the cube case except of the definition of $\beta_{st}$:
\[\cos \frac{s\pi}{l+1}\cos \frac{t\pi}{m+1} +(\cos \frac{s\pi}{l+1}+\cos \frac{t\pi}{m+1})\cosh \beta_{st}=3\]
\subsection{Three dimensional diagonal lattice}
The solution in a 3-dimensional lattice can be obtained by taking $l,m,n,a,b,c,p,q,r \rightarrow \infty$ in the block solution with $ p-a=u, q-b=v, r-c=w $ finite:
\[F_{(0,0,0)}(u,v,w)=\frac{1}{\pi^2}\int_{0}^{\pi} \int_{0}^{\pi} \frac{\cos u \lambda \cos v\mu \exp{(- \left|w \right|\theta)}}{\tanh \theta}d\lambda d\mu \]
where in the cube model we have:
 \[cos\lambda \cos\mu \cosh\theta=1 \]
and in the dodecahedron model we have:
 \[cos\lambda \cos\mu + (\cos\lambda +\cos\mu)\cosh\theta=3 \]
\subsection{Probability of return in 3-dimensional diagonal lattice}
A well known result in case of simple random walk in 3 dimensions is that the probability of return to the starting point is approximately $0.34$; see e.g. McCrea and Whipple [3].\\
We first focus on the probability of return in the diagonal \textit{cube} case.
The probability is $1-\frac{1}{F}$ where
 \[
F_{(0,0,0)}(0,0,0)=\frac{1}{\pi^2}\int_{0}^{\pi} \int_{0}^{\pi} \frac{1}{\tanh \theta}d\lambda d\mu
\]
where  \[\cos\lambda \cos\mu \cosh\theta=1\] \\
Using numerical integration, we find  \[
F_{(0,0,0)}(0,0,0)=\frac{1}{\pi^2}\int_{0}^{\pi} \int_{0}^{\pi} (1-\cos^{2}\lambda \cos^{2}\mu)^{-0.5}d\lambda d\mu \approx 1.3932
\]
In the diagonal cube case we have probability of return $1-\frac{1}{F} \approx 0.2822$ \\

Montroll~\cite{MO} uses a different approach. He observes that many crystals appear as body centered
lattices. The body centered lattice is composed of two interpenetrating simple cubic lattices with the
points of one lattice being at the center of the cubes of the other lattice. The walker moves to one of its
eight neighboring lattice points in the other lattice. He finds the probability of return to the starting point in the diagonal case:$1-\frac{1}{u}\approx.282229985$ where 
$ u=\frac{1}{\pi^3}\int_0^\pi \int_0^\pi
\int_0^\pi(1-\cos{x}  \cos{y} \cos{z})^{-1} dxdydz \approx 1.3932039297 $
\\
We  now focus on the probability of return in the diagonal \textit{dodecahedron} case.

The probability is $1-\frac{1}{F}$ where
 \[
F_{(0,0,0)}(0,0,0)=\frac{1}{\pi^2}\int_{0}^{\pi} \int_{0}^{\pi} \frac{1}{\tanh \theta}d\lambda d\mu
\]
and \[\cos\lambda \cos\mu + (\cos\lambda +\cos\mu)\cosh\theta=3\] \\
Using numerical integration, we find 
 \[
F_{(0,0,0)}(0,0,0)=\frac{1}{\pi^2}\int_{0}^{\pi} \int_{0}^{\pi} [1-(\frac{\cos\lambda +\cos\mu}{3-\cos\lambda \cos\mu})^2]^{-0.5} d\lambda d\mu \approx 1.2298
\]
In the diagonal dodecahedron case we have probability of return $1-\frac{1}{F} \approx 0.1868$ \\

\section{Transformations in two dimensions}
We can transform the diagonal random walk to a simple one by first shrinking with factor $\frac{1}{\sqrt{2}}$ and then a rotation around the origin with angle $\pi /4$.
We get: $(p,q)\rightarrow (\frac{p}{\sqrt{2}},\frac{q}{\sqrt{2}}) \rightarrow (\frac{p-q}{2},\frac{p+q}{2}).$
Let $(x,y)$ be our new coordinate system, then we have: $p=y+x, q=y-x$. We use this transformation to get the desired simple random walk, but now with rotated boundaries.

\subsection{Transformed Rectangular Region}
$I=\{(x,y)| 1\leq y+x\leq m, 1\leq y-x\leq n\}$.\\
 Our original starting point $(a,b)$ is transformed in $(\frac{a-b}{2},\frac{a+b}{2})$.\\
 When starting in $(\frac{a-b}{2},\frac{a+b}{2})$ we get 
 \[
F_1(x,y)= \frac{4}{m+1}\sum_{r=1}^{m} \frac{ \sin\frac{ar\pi}{m+1} \sin\frac{(y+x)r\pi}{m+1} \sinh [(y-x)\beta_r]\sinh[(n+1-b)\beta_r]}{\tanh \beta_r\sinh[(n+1)\beta_r]} 
\]
where $y-x\leq b$.
We prefer to start in $(a,b)$ and then we get:
\begin{multline*}
F_1(x,y)=\\ \frac{4}{m+1}\sum_{r=1}^{m} \frac{ \sin\frac{(a+b)r\pi}{m+1} \sin\frac{(y+x)r\pi}{m+1} \sinh [(y-x)\beta_r]\sinh[(n+1+a-b)\beta_r]}{\tanh \beta_r\sinh[(n+1)\beta_r]} 
\end{multline*}
 \begin{multline*}
F_2(x,y)= \\ \frac{4}{m+1}\sum_{r=1}^{m} \frac{ \sin\frac{(a+b)r\pi}{m+1} \sin\frac{(y+x)r\pi}{m+1} \sinh {[(b-a)\beta_r]} \sinh[(n+1+x-y)\beta_r]}{\tanh \beta_r\sinh[(n+1)\beta_r]} 
\end{multline*}
where $F_1$ is valid for $y-x \leq b-a$ and $F_2$ is valid for $y-x \geq b-a$
and
\[\cos \frac{r\pi}{m+1}\cosh \beta_r=1\]
\subsection{Transformed Semi infinite strip}
$I=\{(x,y)| 1\leq y+x\leq m, 1\leq y-x\}$; we start in $(a,b)$.

\[
F_1(x,y)= \frac{4}{m+1}\sum_{r=1}^{m} \frac{ \sin\frac{(a+b)r\pi}{m+1} \sin\frac{(y+x)r\pi}{m+1} \sinh [(y-x)\beta_r] \exp[(a-b) \beta_r]}{\tanh \beta_r} 
\]
\[
F_2(x,y)= \frac{4}{m+1}\sum_{r=1}^{m} \frac{ \sin\frac{(a+b)r\pi}{m+1} \sin\frac{(y+x)r\pi}{m+1} \sinh [(b-a)\beta_r] \exp[(x-y) \beta_r]}{\tanh \beta_r} 
\]
where $F_1$ is valid for $y-x \leq b-a$ and $F_2$ is valid for $y-x \geq b-a$
and
\[\cos \frac{r\pi}{m+1}\cosh \beta_r=1\]

\subsection{Transformed Infinite strip}
Rotating and shrinking the solution of the semi infinite strip gives, when starting in $(a,a)$:
\[
F(p,s)= \frac{2}{m+1}\sum_{r=1}^{m} \frac{ \sin\frac{2ar\pi}{m+1} \sin\frac{(p+s)r\pi}{m+1} \exp(-\lvert s \lvert\beta_r)}{\tanh \beta_r} \quad (1\leq p+s\leq m)
\]
where \[\cos \frac{r\pi}{m+1}\cosh \beta_r=1\]

\subsection{Transformed Infinite Quadrant}
$I=\{(x,y)| 1\leq y+x, 1\leq y-x\}$; we start in $(a,b)$.
\[
F_1(x,y)= \frac{8}{\pi}\int_{0}^{\pi} \frac{ \sin{[(a+b) \lambda]} \sin{[(y+x) \lambda]} \sinh{[(y-x) \mu]} \exp[(a-b) \mu]}{\tanh \mu} d\lambda 
\]
\[
F_2(x,y)= \frac{8}{\pi}\int_{0}^{\pi} \frac{ \sin{[(a+b) \lambda]} \sin{[(y+x) \lambda]} \sinh{[(b-a) \mu]} \exp[(x-y) \mu]}{\tanh \mu} d\lambda 
\]
where 
$F_1$ is valid for $y-x\leq b-a$ and $F_2$ is valid for $y-x\geq b-a$ and
\begin{equation*}
\cos {\lambda}\cosh {\mu}=1
\end{equation*}
\subsection{Transformed Half-plane}
By taking $m  \rightarrow \infty$ in the solution of the infinite strip, we get when starting in $(a,a)$:
\begin{equation*}
F(p,s)= \frac{2}{\pi}\int_{0}^{\pi} \frac{\sin{(2a \lambda)} \sin{[(p+s) \lambda]} \exp(-\lvert s \lvert\mu)}{\tanh \mu} d\lambda \quad (1\leq p+s)
\end{equation*}
where 
\begin{equation*}
\cos {\lambda}\cosh {\mu}=1
\end{equation*}

\end{document}